\documentclass[12pt]{article}
\usepackage{amsmath, amssymb}
\usepackage{amssymb,pb-diagram}
\usepackage{amscd}

\usepackage{fancybox}

\newcommand{\ZZ}{\mathbb Z}
\newcommand{\PP}{\mathbb P}

\newcommand{\mcT}{\mathcal T}
\newcommand{\mcL}{\mathcal L}

\newcommand{\mcB}{\mathop {\mathcal B}\nolimits}

\newtheorem{thm}{Theorem}[section]
\newtheorem{mthm}{Theorem}

\newtheorem{prop}{Proposition}[section]
\newtheorem{lem}{Lemma}[section]
\newtheorem{defin}{Definition}[section]
\newtheorem{exmple}{Example}[section]
\newtheorem{rem}{Remark}[section]
\newtheorem{qz}{Question}[section]

\renewcommand{\thesubparagraph}{\theparagraph.\@arabic\c@subparagraph}
\begin{document}
%
%
%
%
%

\title{Zariski $N$-ples for a smooth cubic and its tangent lines} 
\author{Shinzo BANNAI and 
Hiro-o TOKUNAGA} 


\maketitle

\begin{abstract}
In this paper, we study the geometry of two-torsion points of elliptic curves in order to distinguish the embedded topology of reducible plane curves consisting of a smooth cubic and its tangent lines. As a result, we obtain a new family of Zariski $N$-ples consisting of such curves.
\end{abstract}

\section*{Introduction}

Let $(\mcB_1, \mcB_2)$ be a pair of reduced complex plane curves in $\PP^2$.  $(\mcB_1, \mcB_2)$ is
said to be a Zariski pair if it satisfies the following two conditons:

\begin{enumerate}
\item[(i)] both $\mcB_1$ and $\mcB_2$ have the same combinatorics 

(see \cite{act} for the details of the
combinatorics of curves),

\item[(ii)] $(\PP^2, \mcB_1)$ is not homeomorphic to $(\PP^2, \mcB_2)$ as a pair of topological spaces.

\end{enumerate}

A $N$-ple of reduced plane curves $(\mcB_1, \ldots, \mcB_N)$ is said to be a Zariski $N$-ple
if  $(\mcB_i, \mcB_j)$ $(1 \le i < j \le N)$ is a Zariski pair.  

The first example of a Zariski pair is given by Zariski in \cite{zariski29}, and for these 25 years much
progress has been made. For example, see \cite{act}. Also see \cite{bgst,  bannai-tokunaga,
bannai-tokunaga17, bty17, bty18,  tokunaga14} for recent results on Zariski pairs for the arrangements consisting of curves of low degree. In particular, in \cite{artal94, bgst, bty18, shirane17},  Zariski pairs for a smooth curves and its tangent lines have been 
studied. In this article, we continue to study such objects: Zariski pairs for a smooth cubic and its tangents.

Note that there exists no Zariski pair for a smooth conic and its $n$ tangent lines as such curves are parametrized by an open set in $\PP^{n-1}$, the set
of  effective divisors of degree $n$. Hence a smooth cubic and its tangent lines is the first object to be studied. In fact, in \cite{artal94, bgst}, Zariski pairs for
an smooth cubic and its inflectional tangent lines are studied. In this note, we study Zariski pairs for a smooth cubic and its $2n$ simple tangents. Let us
explain the combinatorics considered in this article as follows:

Choose distinct points $P_1, \ldots, P_n$ on a smooth cubic, none of which is an inflection point.  For each $P_i$, there exist $4$ lines $L_{P_i, j}$ ($j = 1, 2, 3, 4$) through $P_i$
 which are tangent to $E$ at $Q_{i, j}$ ($j = 1, 2, 3, 4$), respectively.  Choose two of them,  $L_{P_i, j}, L_{P_i, l}$ and put 
 \[
 \mcL_{P_i}^{(j_i, k_i)} = L_{P_i, j_i} + L_{P_i, k_i} \quad \mbox{and} \quad  \mcL = \sum_{i=1}^n\mcL_{P_i}^{(j_i, k_i)}.
 \]
 
 The combinatorics considered in this article is the one given by $E + \mcL$ such that no three lines are concurrent. Now our result can be stated as follows:

\begin{mthm}\label{thm:main} For the combinatorics as above, there exists a Zariski $y(n)$-ple. Here $y(n)$ is the number of $3$-partitions for $n$, i.e.,
it is given as follows:
\[
y(n) = \left \{ \begin{array}{ll}
                         \frac 1{12}(n^2 + 6n + 12) & n \equiv 0 \bmod 6 \\
                         \frac 1{12}(n +1)(n +  5) & n \equiv \pm 1 \bmod 6 \\
                         \frac 1{12}(n + 2)(n + 4 ) & n \equiv \pm 2 \bmod 6 \\
                         \frac 1{12}(n+ 3)^2      & n \equiv 3 \bmod 6
                         \end{array}
                         \right .
\]
\end{mthm}

In previous articles \cite{artal94, bgst},
inflection points which are regarded as three-torsions play  key roles. 
On the other hand, in our proof, a description of torsion points of order $2$ on $E$ plays an important role. More precisely we represent a two-torsion via
intersection points $L_{P_i, j_i}\cap E$,  $L_{P_i, k_i}\cap E$ given in Section 1. This is the new feature in this article.

{\bf Akcnowledgements:}
The first and second authors are  partially supported by Grant-in-Aid for Scientific Research C (18K03263  and 17K05205), respectively.
%
%
%
%
%
%
%
%
%
%
%
%
%
%
%

\section{Preliminaries}

\subsection{Splitting numbers}

In  \cite{shirane2016}, T. Shirane introduced the notion of {\it splitting numbers} and used it to distinguish the embedded topology of curves. In this subsection, we restate the definition and propositions concerning splitting numbers to fit our setting and to simplify the presentation.

Let $B\subset \PP^2$ be a plane curve of even degree and $\phi: X\rightarrow \PP^2$ be the double cover of $\PP^2$ branched along $B$. Let $C\subset\PP^2$ be an irreducible curve.
\begin{defin}
The number of irreducible components of $\phi^\ast(C)$ is called the splitting number of $C$ with respect to $\phi$ and will be denoted by $s_\phi(C)$.
\end{defin}
Note that since we are considering double covers only, $s_\phi(C)=1$ or $2$.
Let $B_1, B_2\subset\PP^2$ be two plane curves of even degree. 

The following proposition allows us to distinguish the embedded topology of curves. For $i=1,2$ let $\phi_i :X_i\rightarrow  \PP^2$ be the double cover branched along $B_i$. Furthermore, let $C_1, C_2\subset \PP^2$ be irreducible curves.

\begin{prop}\label{prop:shirane}
If there exists a homeomorphism $h:\PP^2\rightarrow\PP^2$ with $h(B_1)=B_2$, $h(C_1)=C_2$ then $s_{\phi_1}(C_1)=s_{\phi_2}(C_2)$.
\end{prop}

\begin{proof}
The statement follows directly from \cite[Corollary 1.4]{shirane2016}.
\end{proof}

\subsection{Pairs of tangents and two-torsion points of $E$}\label{subsec:torsion}

Let $E$ be a smooth cubic curve and choose an inflection point $O\in E$.  It is well-known that we can endow $E$ with an abelian group structure on $E$ with $O$ being the zero element  (see
\cite{cassels, silverman-tate}, for example). We denote the addition and subtraction on $E$ by $\dot+$ and $\dot-$. Since we have chosen an inflection point $O$ as the zero element,  for three points $P, Q, R\in E$, $P\dot+Q\dot+ R=O$ if they are collinear. Let $\mcT=\{T_1, T_2, T_3\}$ be the set of non-trivial two-torsion points of $E$. For a point  $P\in E$ which is not an inflection point, it is known that there exists four lines that pass through $P$ and is tangent to $E$ at a point distinct from $P$.  Let,  $L_{P, i}$, $(i=1,\ldots, 4)$ be such four lines and let $Q_i$ be the tangent points. By the geometric description of the group law of elliptic curves, we have
\[
P\dot+2Q_i=O.
\]
Then, for $\{i,j\}\subset\{1,2,3,4\}$ we have
\[
2(P\dot+Q_i\dot+Q_j)=(P\dot+2Q_i)\dot+(P\dot+2Q_j)=O,
\]
hence $P\dot+Q_i\dot+Q_j\in \mcT$. Note that $P, Q_i, Q_j$ cannot be collinear.

\begin{defin}
For a pair  $\mcL_P^{(i,j)}=L_{P,i}+L_{P,j}$ of tangent lines through $P$, the two-torsion point $T= P\dot+Q_i\dot+Q_j$ is called the two-torsion point associated to $\mcL_P^{(i,j)}$.
\end{defin}

\begin{lem}\label{lem:relations}
Under the above setting, 
\begin{align*}
P\dot+Q_1\dot+Q_2=P\dot+Q_3\dot+Q_4=T_1\\
P\dot+Q_1\dot+Q_3=P\dot+Q_2\dot+Q_4=T_2\\
P\dot+Q_1\dot+Q_4=P\dot+Q_2\dot+Q_3=T_3
\end{align*}
for a suitable choice of labels for $Q_i$, $(i=1,2,3, 4)$. Moreover, every non-trivial two-torsion point $T_i$ of $E$ can be obtained as an associated two-torsion point of $\mcL_P^{(i,j)}$ for a suitable choice of pairs of tangent lines.
\end{lem}

\begin{proof}
For $\{i,j\}\subset\{1,2,3,4\}$ we have
\[
2(Q_i\dot-Q_j)=(P\dot+2Q_i)\dot-(P\dot+2Q_j)=O,
\]
hence $Q_i\dot-Q_j$ also becomes a two-torsion point of $E$.

Since $Q_1\dot-Q_2$, $Q_1\dot-Q_3$, $ Q_1\dot-Q_4$  are distinct non-trivial two-torsion points, we can assume that
\[
Q_1\dot-Q_2=T_1, 
Q_1\dot-Q_3=T_2, 
Q_1\dot-Q_4=T_3 
\]
for a suitable choice of labels for  $Q_i$, $(i=1,2,3, 4)$. 
Also, since the subgroup of two-torsion points is isomorphic to $(\ZZ/\ZZ)^{\oplus2}$, we have $T_i\dot+T_j=T_k$ for $\{i,j,k\}=\{1,2,3\}$. These combined with $P\dot+2Q_i=0$ give the desired equalities.
\end{proof}

\section{The case of four tangent lines}

In this section we consider the fundamental case of a smooth cubic and four of its tangent lines. 

Let $P_1, P_2\in E$, $P_1\not=P_2$ be non-inflection points. Then for each $P_i$ ($i=1, 2$),  there exist four lines $L_{P_i,j}$ $(j=1,2,3,4)$ passing through $P_i$ and tangent to $E$ at $Q_{i,j}$ as in Section \ref{subsec:torsion}. We assume that the points  $Q_{i,j}$ are labeled so that at each point $P_i$ the equalities in Lemma \ref{lem:relations} are satisfied.

Let $\mcL^{(i,j), (k,l)}=\mcL_{P_1}^{(i,j)}+\mcL_{P_2}^{(k,l)}=L_{P_1,i}+L_{P_1,j}+L_{P_2,k}+L_{P_2,l}$. Furthermore, let $\phi^{(i,j),(k,l)}: S\rightarrow \PP^2$ be the double cover of $\PP^2$ branched along $\mcL^{(i,j), (k,l)}$. Then we have the following lemma:

\begin{lem}\label{lem:splittingnumber}
Let $T, T^\prime$ be the two-torsion points associated to $\mcL_{P_1}^{(i,j)}, \mcL_{P_2}^{(k,l)}$ respectively, and let $s$ be the splitting number of $E$ with respect to $\phi^{(i,j), (k,l)}$. Then $s=2$ if and only if $T=T^\prime$.
\end{lem}


\begin{proof}
The statement follows from \cite[Proposition 2.5]{shirane2016}. In our case, the divisor $D^\prime_{B,C}$ in  \cite[Proposition 2.5]{shirane2016} coincides with $T+T^\prime$. Then $T=T^\prime$ if and only if the order of $[\mathcal{O}_{\hat{E}} (D^\prime_{B,C})]$ equals 1, which is equivalent to $s=2$ by  \cite[Proposition 2.5]{shirane2016}.
\end{proof}

Now, we consider two curves  $\mcB^1=E+\mcL^{(i_1,j_1), (k_1,l_1)}$ and $\mcB^2= E+\mcL^{(i_2,j_2), (k_2,l_2)}$, each consisting of $E$ and four tangent lines. Note that since $P_1, P_2$ are non-inflectional points, the combinatorics of $\mcB_1$ and $\mcB_2$ are the same for any choice of $(i_1,j_1)$, $(k_1,l_1)$, $(i_2, k_2)$, $(j_2, k_2)\subset\{1,2,3,4\}$. For $\mcB_1, \mcB_2$ we have the following proposition:
\begin{prop}
The pair $(\mcB^1, \mcB^2)$ is a Zariski-pair if the parity of $\left|\{i_1,j_1\}\cap\{k_1,l_1\}\right|$ and  $\left|\{i_2,j_2\}\cap\{k_2,l_2\}\right|$ are distinct, assuming that the labeling of $L_{P_i,j}$ satisfy the equalities in Lemma \ref{lem:relations}. 
\end{prop}
 
 \begin{proof}
Let, $T, T^\prime$ be the  torsion sections associated to $\mcL_{P_1}^{(i,j)}, \mcL_{P_2}^{(k,l)}$ respectively. Then by the equalities in Lemma \ref{lem:relations}, $T=T^\prime$ if and only if $\left|\{i,j\}\cap\{k,l\}\right|$ is even. Now, if the parity of $\left|\{i_1,j_1\}\cap\{k_1,l_1\}\right|$ and  $\left|\{i_2,j_2\}\cap\{k_2,l_2\}\right|$ are distinct, this implies that the splitting number of $E$ with respect to $\phi^{(i_1,j_1), (k_,l_1)}$ and $\phi^{(i_2,j_2), (k_2,l_2)}$ are distinct by Lemma \ref{lem:splittingnumber}. Hence, the pair $(\mcB_1, \mcB_2)$ is a Zariski-pair by \cite[Corollary 1.4]{shirane2016}.
 \end{proof}

\section{Proof of Main Theorem}

For curves of the form $E+\mcL$, we note that if we choose $P_1, \ldots, P_n \in E$ generally, any three of the lines  $L_{P_i,j}$ ($i=1,\ldots, n$, $j=1,\ldots 4$) will not be concurrent.

%
%
Let $\mcT$ be as before and let $\underline{{\rm Sub}}_\wedge(E, \mcL)$ be a set of subarrangements given by
\[
\underline{{\rm Sub}}_\wedge(E, \mcL)=\{E+\mcL_{P_i}^{(k_i,l_i)}\mid i=1, \ldots, n\}.
\]
Define a map $\Phi_\mcL:\underline{{\rm Sub}}_\wedge(E, \mcL) \rightarrow \mcT$ by setting $\Phi_\mcL\left(E+\mcL_{P_i}^{(k_i,l_i)}\right)$ to be the two-torsion associated to $\mcL_{P_i}^{(k_i,l_i)}$. 
With $\Phi_\mcL$, we have a $3$-partition of $\underline{{\rm Sub}}_\wedge(E, \mcL)$ by $\bigcup_i \Phi_\mcL^{-1}(T_i)$. In the following we denote the subarrangement $E+\mcL_{P_i}^{(k_i,l_i)}$ by $[P_i, k_i, l_i]$ to simplify the notation.

\begin{defin}
Under the above settings, the 3-partition $(m_1, m_2, m_3)$ of $n$ associated to $E+\mcL$ is defined to be a triple of non-negative integers $(m_1, m_2, m_3)$ such that $m_1\geq m_2 \geq m_3$ and $\{m_1, m_2, m_3\}=\left\{\left|\Phi_\mcL^{-1}(T_1)\right|, \left|\Phi_\mcL^{-1}(T_2)\right|, \left|\Phi_\mcL^{-1}(T_3)\right|\right\}$. 
\end{defin}

Note that in the above definition, $m_1+m_2+m_3=n$. The integer $m_i$ need not be equal to $\left|\Phi_\mcL^{-1}(T_i)\right|$, the labels may be rearranged. Also, Lemma \ref{lem:relations} implies that every 3-partition $(m_1, m_2, m_3)$ of $n$ can be obtained as a 3-partition associated to $E+\mcL$ by choosing $\mcL_{P_i}^{(k_i,l_i)}$ suitably.

Now, Theorem \ref{thm:main} follows from the following proposition. 

\begin{prop}
Let $\mcB^1=E+\mcL$, $ \mcB^2=E+\mcL^\prime$ and $(m_1, m_2, m_3)$, $(m^\prime_1, m^\prime_2, m^\prime_3)$ be the associated 3-partitions of $n$ respectively. If there exists a homeomorphism $h:\PP^2\rightarrow\PP^2$ such that $h(\mcB^1)=\mcB^2$, then $(m_1, m_2, m_3)=(m^\prime_1, m^\prime_2, m^\prime_3)$.
\end{prop}

\begin{proof}
Suppose there exists an homeomorphism $h:\PP^2\rightarrow \PP^2$ such that $h(\mcB^1)=\mcB^2$. Then $h$ naturally induces a bijection $h_\natural:   \underline{{\rm Sub}}_\wedge(E, \mcL_1) \rightarrow   \underline{{\rm Sub}}_\wedge(E, \mcL_2)$. Furthermore, $\Phi_\mcL([P_{i_1}, k_{i_1}, l_{i_1}])=\Phi_\mcL([P_{i_2}, k_{i_2}, l_{i_2}])$ if and only if $\Phi_{\mcL^\prime}(h_\natural([P_{i_1}, k_{i_1}, l_{i_1}]))=\Phi_{\mcL^\prime}(h_\natural([P_{i_2}, k_{i_2}, l_{i_2}]))$ since the splitting number of $E$ with respect to the double cover branched along  
\[L_{P_{i_1},k_{i_1}}+L_{P_{i_1},l_{i_1}}+L_{P_{i_2},k_{i_2}}+L_{P_{i_2},l_{i_2}} \text{ and } h(L_{P_{i_1},k_{i_1}}+L_{P_{i_1},l_{i_1}}+L_{P_{i_2},k_{i_2}}+L_{P_{i_2},l_{i_2}})\]
must be equal by Proposition \ref{prop:shirane} and Lemma \ref{lem:splittingnumber}. Moreover, $h$ naturally induces a bijection  $h_\flat:\mcT \rightarrow \mcT$ such that the following diagram commutes.
\[
\begin{CD}
 \underline{{\rm Sub}}_\wedge(E, \mcL_1) @>\Phi_{\mcL_1}>> \mcT \\
 @V h_\natural VV  @VVh_\flat V \\
 \underline{{\rm Sub}}_\wedge(E, \mcL_1) @> \Phi_{\mcL_2} >> \mcT
\end{CD}
\]
Hence, we have $(m_1, m_2, m_3)=(m^\prime_1, m^\prime_2, m^\prime_3)$.
\end{proof}

\noindent Shinzo BANNAI\\
National Institute of Technology, Ibaraki College\\
866 Nakane, Hitachinaka-shi, Ibaraki-Ken 312-8508 JAPAN \\
{\tt sbannai@ge.ibaraki-ct.ac.jp}\\

\noindent Hiro-o TOKUNAGA\\
Department of Mathematical  Sciences, Graduate School of Science, \\
Tokyo Metropolitan University, 1-1 Minami-Ohsawa, Hachiohji 192-0397 JAPAN \\
{\tt tokunaga@tmu.ac.jp}

\end{document}